\documentclass{amsart}

\usepackage{amsmath,amsthm,amsopn,amstext,amscd,amsfonts,amssymb,mathrsfs,mathtools}
\usepackage{dsfont}
\usepackage{comment}
\usepackage{xcolor}
\usepackage{float}
\usepackage[active]{srcltx}
\usepackage{graphicx, mathdots}

  \usepackage{mathtools}

\newtheorem*{theorem}{\sc Theorem}
\newtheorem*{conj}{\sc Conjecture}

\newtheorem*{eje}{\sc Example}
\newtheorem*{coro}{\sc Corollary}

\usepackage[active]{srcltx}

\newcommand{\dps}{\displaystyle}
\usepackage{graphicx, epsfig, subfig}
\usepackage{lscape}
\usepackage{rotating}
\usepackage{caption}
\usepackage{multicol}
\usepackage{colortbl}
\usepackage{nicematrix}

\newmuskip\pFqskip
\mathchardef\pFcomma=\mathcode`, 

\makeatletter
\def\BState{\State\hskip-\ALG@thistlm}
\makeatother

\def\downbar#1{
\setbox10=\hbox{$#1$}
            \dimen10=\ht10 \advance\dimen10 by 2.5pt
            \ifdim \dimen10<15pt 
               \advance\dimen10 by -0.5pt
               \dimen11=\dimen10
               \advance\dimen10 by 2.5pt
               \lower \dimen11
            \else \lower \ht10 \fi
            \hbox {\hskip 1.5pt \vrule height \dimen10 depth \dp10}}
\def\upbar#1{
\setbox10=\hbox{$#1$}
            \dimen10=\ht10 \advance\dimen10 by \dp10 \advance\dimen10 by 2.5pt
            \ifdim \dimen10<15pt 
                \advance\dimen10 by 2pt \fi
            \raise 2.5pt \hbox {\hskip -1.5pt \vrule height \dimen10}}

\begin{document}
\title[On the positivity of a certain function]{On the positivity of a certain function related with the Digamma function}
\author{K. Castillo}
\address{CMUC, Department of Mathematics, University of Coimbra, 3001-501 Coimbra, Portugal}
\email{ kenier@mat.uc.pt}

\date{\today}
\begin{abstract}
It is proved that
$$
\left(\frac{x^n}{1-e^{-x}}\right)^{(n)}>0
$$
for all $x\in (\log 2, \infty)$ and $n\in \mathbb{N}$,  which improves the result of [Al-Musallam and Bustoz in Ramanujan J. 11 (2006) 399-402].
\end{abstract}
\maketitle
\section{Introduction} 
A function $f: (a, b)\subset \mathbb{R}\longrightarrow \mathbb{R} $ is completely monotonic if it is infinitely differentiable and 
$$
(-1)^n f^{(n)}(x)\geq 0
$$
for all $x\in (a, b)$ and $n \in \mathbb{N}$. A function $f(-x)$ is called absolutely monotonic on $(-b, -a)$ if and only if $f(x)$ is completely monotonic on $(a, b)$.  Absolutely monotonic functions were pioneeringly introduced by Bernstein. Bernstein himself, and later Widder independently, discovered that a necessary and sufficient condition for $f$ to be completely monotonic on $(0, \infty)$ is that 
$$
f(x)=\mathcal{L}(\mu)(x)=\int e^{-xt}\,\mathrm{d}\mu(t),
$$
where $\mu$ is a positive measure on $[0, \infty)$ and the integral converges for all positive $x$. (These and other classical results on absolutely/completely monotonic functions can be found in \cite[Chapter IV]{W46} and \cite{B71}.) As it was remarked in \cite{ABK05}, by Bernstein's theorem, its is easy to see that the absolute value of the digamma function, $\psi=\Gamma'/\Gamma$, and the absolute value of  its derivatives (the polygamma functions) are completely monotonic functions on $(0, \infty)$. Indeed, 
\begin{align*}
(-1)^{n+1}\psi^{(n)}(x)=\mathcal{L} \left( \frac{t^n}{1-e^{-t}} \right)(x)
\end{align*}
 for all $x\in (0, \infty)$ and $n\in \mathbb{N}$.

In \cite{CI03} Clark and Ismail introduced the functions
$$
F_m(x)=x^m \psi (x), \quad G_m(x)=-x^m \psi(x).
$$
They proved that $F_m^{(m+1)}$ is completely monotonic on $(0, \infty)$ for $m\in\mathbb{N} \setminus\{0\}$ \cite[Theorem 1.2]{CI03}  and that $G_m^{(m)}$ is completely monotonic on $(0, \infty)$ for $m=1, 2,\dots, 16$  \cite[Theorem 1.3]{CI03}. Afterwards, they wrote: {\it ``We believe Theorem $1.3$ $[$$G_m^{(m)}$ is completely monotonic on $(0, \infty)$$]$ is true for all $m$ $[$...$]$".} However, Alzer, Berg, and Koumandos \cite[Theorem 1.1]{ABK05} proved that there exists an integer $m_0$ such that for all $m\geq m_0$ the function  $G_m^{(m)}$ is not completely monotonic.  From this and the relation \cite[(2.4)]{CI03}
\begin{align}\label{aux1}
G^{(m)}_m(x)= \mathcal{L}\left(t^m \left(\frac{t^m}{1-e^{-t}}\right)^{(m)}\right)(x),
\end{align}
it follows  that the following conjecture of Clark and Ismail \cite[Conjecture 1.4]{CI03} is false:

\begin{conj}
\begin{align}\label{pos}
\left(\frac{x^n}{1-e^{-x}}\right)^{(n)}>0
\end{align}
for all $x\in (0, \infty)$ and $n\in \mathbb{N}$.
\end{conj}
By showing that \eqref{pos} holds for $n=1, 2,\dots, 16$ (and using \eqref{aux1}), Clark and Ismail proved that $G_n^{(n)}$ is (strictly) completely monotonic on $(0, \infty)$ for these values of $n$. Regardless of the fact that the conjecture is not true, the inequality \eqref{pos} is of interest in its own right. It remains an open problem to determine the smallest positive number $a$ (positive integer $n_0$) such that \eqref{pos} remains positive for all $x\in (a, \infty)$ and $n\in \mathbb{N}$ ($x\in (0, \infty)$ and $n>n_0$ with $n\in \mathbb{N}$). (This open problem was also placed in \cite[Section 4]{ABK05}.) In \cite[Theorem 2.1]{AB03}\footnote{This paper was submitted on May 21, 2003 and accepted for publication on October 24, 2003. J. Bustoz passed away on August 13, 2003.}, Al-Musallam and Bustoz proved that \eqref{pos} holds for all $x\in (2 \log 2, \infty)$ and $n\in \mathbb{N}$. (This was also proved independently in \cite[p. 112]{ABK05} using the same idea: an inequality proved by Szeg\H{o}  \cite[Theorem 17a, p. 168]{W46}.)  Our main theorem, which improves the result in \cite{AB03}, reads as follows:

\begin{theorem}
 \eqref{pos} holds for all $x\in (\log 2, \infty)$ and $n\in \mathbb{N}$. 
\end{theorem}


As in \cite[Theorem 3.1]{AB03}, now using the above theorem, the next result follows. (The details are left to the reader.)

\begin{coro}
\begin{align*}
\left(\frac{x^{n+\alpha}}{1-e^{-x}}\right)^{(n)}>0
\end{align*}
for all $\alpha \in (0, \infty)$,  $x\in (\log 2, \infty)$, and $n\in \mathbb{N}$.
\end{coro}

\begin{eje}
It is easy to obtain from \cite[(1), p. 11]{E01}, for  $n\in \mathbb{N}\setminus\{0\}$, the power series
$$
\left(\frac{x^{n}}{1-e^{-x}}\right)^{(n)}=\frac{n!}{2}+\sum_{j=2}^{\infty} \frac{(j+n-1)!}{(j-1)!} \frac{B_j}{j!} x^{j-1}
$$
valid in the disk $|x|<2\pi$ which extends to the nearest singularities $x=\pm 2 \pi i$ of $x/(e^x-1)$. (The coefficients $B_j$ are the Bernoulli numbers. The odd Bernoulli numbers are all zero after the first, but it is a highly complex task to determine the even Bernoulli numbers.) Let us imagine that we are questioned about the sign of the following sum:
\begin{align*}
S_n&=\frac{(n+0)!}{0!\,1!} B_0+\frac{(n+1)!}{1!\,2!} B_2+\frac{(n+3)!}{3!\,4!} B_4+\frac{(n+5)!}{4! \, 5!} B_6+\cdots\\[7pt]
&=n!+\sum_{j=1}^{\infty} \frac{(2j+n-1)!}{(2j-1)!} \frac{B_{2j}}{(2j)!}.
\end{align*}
(Recall that $B_0=1$.) Note that
$$
\left.\left(\frac{x^{n}}{1-e^{-x}}\right)^{(n)}\right|_{x=1}=S_n-\frac{n!}{2}.
$$
Since $\log2\approx 0.693147 <1<2\pi$, our main results gives 
$$
S_n>n!/2
$$ 
 $n\in \mathbb{N}\setminus\{0\}$. It is worth pointing out that from the results obtained in \cite{AB03, ABK05}, it is not possible to conclude this because $2 \log 2 \approx 1.38629>1$. Now it only remains to check that $S_n$ converges, which follows from
$$
\lim_{j\to \infty}\sqrt[\leftroot{-2}\uproot{2}\dps 2j]{\frac{(2j+n-1)!}{(2j-1)!} \frac{|B_{2j}|}{(2j)!}}=\frac{1}{2\pi}<1.
$$
In \cite{ABK05} the relation of the function given in \eqref{pos} with a function of Hardy and Littlewood was extensively explored.
\end{eje}
\section{Proof of the theorem}
Set
$$
f_n(x)=\frac{\mathrm{d}^n}{\mathrm{d} x^n}\left(\frac{x^n}{1-e^{-x}}\right).
$$
If $c>0$ is arbitrary and fixed, the series 
$$
\frac{1}{1-e^{-x}}=\sum_{j=0}^\infty e^{-j x},
$$
converges uniformly on $[c, \infty)$. We then write $f_n$ in the form 
$$
f_n(x)=\sum_{j=0}^\infty \dps \frac{\mathrm{d}^n}{\mathrm{d} x^n}\left( e^{-j x} x^n\right).
$$
Recall that \cite[(5), p. 188]{B81} $n! L_n(x)=e^x (\mathrm{d}^n/\mathrm{d}x^n)\left( e^{-x} x^n\right)$, 
$L_n$ being the Laguerre polynomial of degree $n$, and so
$$
n! e^{-jx}L_n(jx)=\frac{\mathrm{d}^n}{\mathrm{d} x^n}\left( e^{-j x} x^n\right).
$$
Hence 
\begin{align*}
f_n(x)=n! \sum_{j=0}^\infty L_n(jx) e^{-jx}
\end{align*}
on $[c, \infty)$. There is a well-known connection between the Laguerre and Hermite polynomials due to Feldheim \cite[(33), p. 195]{B81}:
$$
\int_0^\infty e^{-t^2} H_n^2(t) \cos(2^{1/2} y\,t)\,\mathrm{d}t=\sqrt{\pi} 2^{n-1} n! L_n(y^2),
$$
$H_n$ being the Hermite polynomial of degree $n$. Write
$$
y^2=j\,x.
$$
From the above expressions, we have
$$
f_n(x)=\frac{1}{\dps \sqrt{\pi}\, 2^{n-1} } \sum_{j=0}^\infty \int_0^\infty g_j(t)\, \mathrm{d}t,
$$
where
$$
g_j(t)=e^{-t^2} e^{-jx} H_n^2(t) \cos (\sqrt{2 j\, x} \,t)
$$
for all $x\in [c, \infty)$. (Recall that it is not true that uniform convergence is sufficient to allow the interchanging of the sum and integral when the integral is over an infinity interval.)
However, the function $g_j$ is integrable and
$$
 \sum_{j=0}^\infty \int_0^\infty|g_j(t)|\,\mathrm{d}t<0,
 $$
for all $x\in [c, \infty)$. These conditions allow the interchanging of the above sum and integral, see, for instance, \cite[Corollary 17.4.7]{T19}. Indeed, since
$$
|g_j(t)|<e^{-t^2} e^{-j\,x} H_n^2(t),
$$
we see at once that 
\begin{align*}
\int_0^\infty|g_j(t)|\,\mathrm{d}t&< e^{-j\,x}\, \int_0^\infty e^{-t^2} H_n^2(t)\,\mathrm{d}t\\[7pt]
&< e^{-j\,x}\, \int_{-\infty}^\infty e^{-t^2} H_n^2(t)\,\mathrm{d}t\leq\sqrt{\pi} 2^n n!\, e^{-j\,c}<\infty,
\end{align*}
and the integrability of $g_j$ is guaranteed. Moreover, 
\begin{align*}
\sum_{j=0}^\infty \int_0^\infty|g_j(t)|\,\mathrm{d}t&<\sqrt{\pi}\, 2^n\, n!\,\sum_{j=0}^\infty e^{-j\,x}\\[7pt]
&\leq\sqrt{\pi}\, 2^n\,n!\,\dps \frac{e^c}{e^{c}-1}<\infty.
\end{align*}
Consequently, we can interchange the sum and integral to obtain
$$
\sqrt{\pi} 2^{n-1} f_n(x)=\int_0^\infty e^{-t^2} H_n^2(t) \sum_{j=0}^\infty e^{-jx} \cos (\sqrt{2 j\, x} \,t)\,\mathrm{d}t.
$$
Finally, note that
$$
 \sum_{j=0}^\infty e^{-jx} \cos (\sqrt{2 j\, x} \,t)>1- \sum_{j=1}^\infty e^{-jx} =1-\frac{1}{e^{x}-1}=g(x).
$$
Thus $g(x)\geq 0$ if and only if $x>\log 2$. This completes the proof.

\section*{Acknowledgements}

This work was supported by the Centre for Mathematics of the University of Coimbra-UIDB/00324/2020, funded by the Portuguese Government through FCT/ MCTES. 


\begin{thebibliography}{10}

\bibitem{AB03} F. Al-Musallam, J. Bustoz, On a conjecture of W. E. Clark and M. E. H. Ismail, Ramanujan J., 11 (2006) 399-402.

\bibitem{ABK05} H. Alzer, C. Berg, S. Koumandos, On a conjecture of Clark and Ismail, J. Approx. Theory, 134 (2005) 102-113.

\bibitem{B71} R. P. Boas, Signs of derivatives and analytic behavior, Amer. Math. Monthly, 78 (1971) 1085-1093.

\bibitem{CI03} W. E. Clark, M.E.H. Ismail, Inequalities involving gamma and psi functions, Anal. Appl., 1 (2003) 129-140.

\bibitem{E01} H. M. Edwards, Riemann's zeta function. Reprint of the 1974 origina. Dover Publications, Inc., Mineola, NY, 2001. 

\bibitem{B81}  A. Erd\'elyi, W. Magnus, F. Oberhettinger, F. G. Tricomi, Higher transcendental functions. Vol. II. Based on notes left by Harry Bateman. Reprint of the 1953 original. Robert E. Krieger Publishing Co., Inc., Melbourne, Fla., 1981. xviii+396 pp. 

 \bibitem{T19} W. J. Terrell, A passage to modern analysis. Pure and Applied Undergraduate Texts, 41. American Mathematical Society, Providence, RI, 2019. xxvii+607 pp. 

\bibitem{W46}  D. V. Widder, The Laplace Transform, Princeton University Press, Princeton, 1946.


\end{thebibliography}
 \end{document}